\documentclass[10pt,a4paper]{amsart}
\usepackage[latin1]{inputenc}
\usepackage[T1]{fontenc}
\usepackage[english]{babel}

\usepackage{amsmath,amssymb,amsthm}

\theoremstyle{plain}
\newtheorem{theorem}{Theorem}[section]

\newtheorem{prop}[theorem]{Proposition}

\theoremstyle{definition}
\newtheorem{remark}[theorem]{Remark}
\newtheorem{example}[theorem]{Example}
\newtheorem{examples}[theorem]{Examples}
\newtheorem{definition}[theorem]{Definition}

\newcommand{\C}{\mathbb{C}}
\newcommand{\R}{\mathbb{R}}
\newcommand{\K}{\mathbb{K}}
\newcommand{\N}{\mathbb{N}}
\newcommand{\T}{\mathbb{T}}
\newcommand{\eps}{\varepsilon}
\newcommand{\ext}[1][X^*]{\ensuremath{\mathrm{ext}(B_{#1})}}
\newcommand{\extr}{\ensuremath{\mathrm{ext}}}
\newcommand{\e}{\mathrm{e}}
\newcommand{\Id}{\mathrm{Id}}

\renewcommand{\leq}{\leqslant}
\renewcommand{\geq}{\geqslant}

 \DeclareMathOperator{\re}{Re}
 \DeclareMathOperator{\dist}{dist}
 \DeclareMathOperator{\ecc}{\overline{\mathrm{co}}}
 \DeclareMathOperator{\ec}{co}

\begin{document}
\title{Numerical index of Banach spaces and duality}
 \date{August 1st, 2005}
 \thanks{The work of the second-named author
was supported by a fellowship from the
\textit{Alexander-von-Humboldt Stiftung}. The work of the
third-named author was partially supported by Spanish MCYT project
no.\ BFM2003-01681 and Junta de Andaluc\'{\i}a grant FQM-185.}

\maketitle

\vspace{0.3cm}

\markboth{Numerical index and duality}{K.~Boyko, V.~Kadets,
M.~Mart\'{\i}n, D.~Werner}

\centerline{\textsc{\large Konstantin Boyko, \quad Vladimir
Kadets,}}

\begin{center}\small Faculty of Mechanics and Mathematics \\ Kharkov National
University\\ pl.\ Svobody 4, 61077 Kharkov, Ukraine \\
\emph{E-mail:} \texttt{k\_boyko@ukr.net, vovalkadets@yahoo.com}
\end{center}

\vspace{0.3cm}

\centerline{\textsc{\large Miguel Mart\'{\i}n\footnote{Corresponding
author.},}}

\begin{center}\small Departamento de An\'{a}lisis Matem\'{a}tico \\ Facultad de
Ciencias \\ Universidad de Granada \\ 18071 Granada, Spain \\
\emph{E-mail:} \texttt{mmartins@ugr.es}
\end{center}

\vspace{0.3cm}

\centerline{\textsc{\large and Dirk Werner}}

\begin{center}\small
Department of Mathematics\\ Freie Universt\"{a}t Berlin\\ Arnimallee
2-6, D-14195 Berlin, Germany \\ \emph{E-mail:}
\texttt{werner@math.fu-berlin.de}
\end{center}

  \thispagestyle{empty}

\begin{abstract}
We present an example of a Banach space whose numerical index is
strictly greater than  the numerical index of its dual, giving a
negative answer to a question which has been latent since the beginning
of the seventies. We also show a particular case in which the
numerical index of the space and the one of its dual coincide.
\end{abstract}

\newpage

\section{Introduction}
The concept of numerical index of a Banach space was first suggested
by G.~Lumer in 1968 (see \cite{D-Mc-P-W});  it is a parameter
relating the norm and the numerical range of operators on the space.
The notion of numerical range was first introduced by O.~Toeplitz in
1918 \cite{Toe} for matrices, and it was extended in the sixties to
bounded linear operators on an arbitrary Banach space by F.~Bauer
\cite{Bauer} and G.~Lumer \cite{Lumer}. Classical references here
are the monographs by F.~Bonsall and J.~Duncan \cite{B-D1,B-D2}. For
recent results we refer the reader to
\cite{Ed-Dari,Ed-Dari-Khamsi,F-M-P,K-M-RP, MarRNP,
MarMer,MaMeRo,M-V,Oik04}, and to the expository paper \cite{Mar} and
references therein.

Here and subsequently, for a real or complex Banach space $X$, we
write $B_X$ for the closed unit ball and $S_X$ for the unit sphere
of $X$. The dual space is denoted by $X^*$, and  the Banach algebra
of all bounded linear operators on $X$ by $L(X)$. The
\emph{numerical range} of such an operator $T$ is the subset $V(T)$
of the scalar field defined by
$$
V(T):=\{x^*(Tx)\ : \ x\in S_X,\ x^*\in S_{X^*},\ x^*(x)=1\}.
$$
The \emph{numerical radius} of $T$ is the seminorm defined on $L(X)$
by
$$
v(T):=\sup\{|\lambda|\ : \ \lambda\in V(T)\}
$$
for each $T\in L(X)$. The \emph{numerical index} of the space $X$,
is the constant $n(X)$ defined by
$$
n(X):=\inf\{v(T)\ : \ T\in L(X),\ \|T\|=1\}
$$
or, equivalently, the greatest constant $k\geq 0$ such that
$k\|T\|\leq v(T)$ for every $T\in L(X)$. Note that $0\leq n(X) \leq
1$, and $n(X)>0$ if and only if $v$ and $\|\cdot\|$ are equivalent
norms on $L(X)$ (the numerical radius can be a non-equivalent norm
on $L(X)$; see \cite[Example~3.b]{M-P}). In the complex case, it
is a celebrated result due to H.~Bohnenblust and S.~Karlin
\cite{B-K} (see also \cite{Gli}) that $n(X) \geq 1/\e$, so the
numerical radius is always an equivalent norm. Actually, the set of
values of the numerical index was established by J.~Duncan,
C.~McGregor, J.~Pryce, and A.~White \cite{D-Mc-P-W}, who proved that
\begin{eqnarray*}
\{n(X)\ : \ X \text{ complex Banach space} \ \}&=&[\e^{-1},1], \\
\{n(X)\ : \ X \text{ real Banach space} \ \}&=&[0,1].
\end{eqnarray*}

Even before the name of numerical index was introduced, it was known
that a Hilbert space of dimension greater than one has numerical
index $1/2$ in the complex case, and $0$ in the real case (see
\cite[\S17]{Hal}). $L$- and $M$-spaces have numerical index~$1$
\cite{D-Mc-P-W}, a property shared by the disk algebra
\cite[Theorem~3.3]{C-D-Mc}, and by every Banach space nicely
embedded into any $C_b(\Omega)$-space \cite[Corollary~2.2]{WerJFA}
(even by every space that is
semi-nicely embedded  into any $C_b(\Omega)$-space
\cite[Corollary~2]{MarRNP}). Very recently, approximations to
the computation of the numerical index of the $L_p(\mu)$-spaces have
been made \cite{Ed-Dari,Ed-Dari-Khamsi}, and the exact computation
of the numerical indices of the two-dimensional spaces whose unit
balls are regular polygons appears in \cite{MarMer}.

Let us mention here a couple of facts concerning the numerical index
which will be relevant to our discussion. Let us fix a bounded
linear operator $T$ on a Banach space $X$. It is a well-known result
of the theory of numerical ranges (see \cite[\S 9]{B-D1}) that
$$
\sup \re V(T)=\lim_{\alpha\downarrow 0} \dfrac{\|\Id  + \alpha\, T\|
-1}{\alpha}
$$
and so,
$$
v(T)=\max_{\omega\in \T}\,\lim_{\alpha\downarrow 0} \dfrac{\|\Id  +
\alpha\,\omega\, T\| -1}{\alpha},
$$
where $\T$ stands for the unit sphere of the base field $\K$ ($=\R$
or $\C$). On the one hand, we can deduce from the above formula that
$$
v(T)=\|T\| \qquad \Longleftrightarrow \qquad \max_{\omega\in \T}
\|\Id  + \omega\,T\|= 1 + \|T\|
$$
(see \cite[Lemma~2.3]{MaOi}). On the other hand, it also implies
that $v(T)=v(T^*)$, where $T^*$ is the adjoint operator of $T$, and
it clearly follows that
\begin{equation}\label{eq:inequality}\tag{$\ast$}
n(X^*)\leq n(X)
\end{equation}
for every Banach space $X$ (see \cite[\S 32]{B-D2}). The question if
this is actually an equality, which is certainly true for reflexive
spaces, has been around  from the beginning of the
subject.

The main aim of this paper is to give a negative answer to the above
question, i.e., we will show that the numerical index of the dual of
a Banach space can be strictly smaller than the numerical index of
the space.

The outline of the paper is as follows. In Section~\ref{sec:C-rich}
we introduce a massiveness property for a Banach space
called ``lushness''  which
implies numerical index~$1$, and we prove that C-rich subspaces of
$C(K)$-spaces satisfy it. Next, we use the above results in
Section~\ref{sec:counterexample} to present examples of Banach
spaces whose numerical index is strictly bigger than the numerical
index of their duals, and other related counterexamples. Finally, we
devote Section~\ref{sec:positive} to show a positive result: the
dual of a Banach space having the Radon-Nikod\'{y}m property (RNP for
short) and numerical index~$1$ also has numerical index~$1$.

We finish this introduction by recalling some definitions and
fixing notation.

Let $X$ be Banach space. Recall that $x_0\in B_X$ is said to be a
\emph{denting point} of $B_X$ if it belongs to slices of $B_X$
with arbitrarily small diameter. More precisely, for each $\eps>0$
one can find a functional $x^*\in S_{X^*}$ and a positive number
$\alpha$ such that the slice
$$
S(B_X,x^*,\alpha):=\{x\in B_X\ : \ \re x^*(x)>1-\alpha\}
$$
contains $x_0$ and is contained in turn in the closed ball centered
at $x_0$ with radius $\eps$. If $X$ is a dual space and the
functionals $x^*$ can be taken to be $w^*$-continuous, then we say
that $x_0$ is a \emph{$w^*$-denting point}. If $B$ is a subset of
$X$, we write $\ec(B)$ and $\ecc(B)$ to denote, respectively, the
convex and closed convex hull of $B$. Then, $\ec(\T\,B)$ will be the
absolutely convex hull of $B$. Finally, we denote by $\extr(A)$ the
set of extreme points of the convex subset $A\subseteq X$.

\section{Lush spaces and C-rich subspaces of $C(K)$}\label{sec:C-rich}

A Banach space $X$ is an \emph{almost-CL-space} if $B_X$ is the
closed absolutely  convex hull of every maximal convex subset of
$S_X$. This notion was introduced by \AA.~Lima \cite{Lima2},
generalizing the concept of \emph{CL-space} (the same definition
without closure) given by R.~Fullerton \cite{Full} in 1960. We refer
to \cite{MartPaya-CL,Reis} and references therein for recent
results. Real and complex almost-CL-spaces have numerical index~$1$
(see \cite[\S4]{Mar}). Actually, the basic examples of Banach spaces
with numerical index~$1$ are known to be almost-CL-spaces (see
\cite{MartPaya-CL} and \cite[Theorem~32.9]{B-D2}). The next
definition is a weakening of the concept of almost-CL-space which
still implies numerical index~$1$. We will show later
(Example~\ref{example:Asplund}) that this weakening is strict,
giving in particular an example of a Banach space with numerical
index~$1$ which is not an almost-CL-space.

\begin{definition} \label{deffat}
We say that a Banach space $X$ is \emph{lush} if for every $x ,y\in
S_X$ and every $\eps>0$, there exists $y^* \in S_{Y^*}$  such that
$y \in S(B_X,y^*, \eps)$ and
$$
\dist\bigl(x,\ec\bigl(\T\, S(B_X,y^*,\eps)\bigr)\bigr) < \eps.
$$
\end{definition}

The (immediate) proof of the fact that almost-CL-spaces have
numerical index~$1$ can be straightforwardly extended to lush
spaces.

\begin{prop}
Let $X$ be a lush Banach space. Then $n(X)=1$.
\end{prop}

\begin{proof}
For $T\in L(X)$ with $\|T\|=1$, and $0<\eps<1/2$ fixed, we take
$x_0\in S_X$ such that $\|Tx_0\|>1-\eps$, and we apply the
definition of lushness to $x_0$ and $y_0=\dfrac{Tx_0}{\|Tx_0\|}$ to
get $y^* \in S_{Y^*}$ with $y_0 \in S(B_X,y^*, \eps)$ and $x_1,
\ldots, x_n \in S(B_X,y^*, \eps)$, $\theta_1, \ldots, \theta_n\in
\T$ such that a convex combination $v=\sum \lambda_k\theta_k x_k$ of
elements $\theta_1 x_1, \ldots, \theta_n x_n$ approximates $x_0$ up
to $\eps.$ Then
$$
|y^*(Tv)| = \left|y^*(y_0) - y^*\left(T\left(\dfrac{x_0}{\|Tx_0\|} -
v\right)\right)\right|
> 1 - 4\eps,
$$
but on the other hand $y^*(Tv)$ is a convex combination of
$y^*(\theta_1T x_1), \ldots, y^*(\theta_n Tx_n).$ So there is an
index $j$ such that
$$
|y^*(Tx_j)| = |y^*(\theta_j Tx_j)|
> 1 - 4\eps.
$$
Now, we have
\begin{align*}
\max_{\omega\in \T} \|\Id  + \omega\,T\|  & \geq \max_{\omega\in \T}
\left|y^*\big([\Id  + \omega\, T](x_j)\big)\right|
 \geq \max_{\omega\in \T} \left|y^*(x_j) + \omega y^*(T x_j)\right|
\\ & = |y^*(x_j)| + |y^*(T x_j)| > 2-5\eps.
\end{align*}
Letting $\eps\downarrow 0$ we deduce that $\displaystyle
\max_{\omega\in \T} \|\Id  + \omega\,T\|= 1 + \|T\|$ and therefore,
$v(T)=\|T\|$.
\end{proof}

Real or complex $C(K)$-spaces are almost-CL-spaces (actually, they
are CL-spaces, see \cite{MartPaya-CL}) and therefore, lush. We now present
a wide class of subspaces of $C(K)$ which are lush, but, as we
will show in Example~\ref{example:Asplund}, they are not
almost-CL-spaces in general.

\begin{definition}
Let $K$ be a compact Hausdorff  space. A closed subspace
$X$ of $C(K)$ is said to be \emph{C-rich} if for every nonempty open
subset $U$ of $K$ and every $\eps > 0$, there is a positive function
$h$ of norm $1$ with support inside $U$ such that the distance from $h$
to $X$ is less than~$\eps$.
\end{definition}

\begin{theorem}\label{th:C-rich=>fat}
Let $K$ be a compact Hausdorff   space and let $X$ be a
C-rich subspace of $C(K)$. Then $X$ is lush and, therefore,
$n(X)=1$.
\end{theorem}

\begin{proof}
We fix $x,y\in S_X$ and $\eps>0$. We take $t_0\in K$ such that
$|y(t_0)|=1$ and we write $a=x(t_0)$, $b = y(t_0)$. Find an open
subset $U$ of $K$ with $t_0\in U$ and such that
\begin{equation}\label{eq:c-rich-u-small}
|x(t)-a|< \eps/4 \qquad \text{and} \qquad |y(t)-b|<\eps/4
\end{equation}
for every $t\in U$. Finally, the C-richness of $X$ gives us a
norm-one function $h:K\longrightarrow [0,1]$ with support inside $U$
and distance to $X$ less than $\eps/4$. Let $\tilde h \in S_X$ be a
function with
\begin{equation}\label{eq:tildeh}
\|\tilde h - h\| < \eps/4.
\end{equation}
Since $\|h\|=1$, there is $t_1\in U$ such that $h(t_1)=1$ and, by
Eq.~\eqref{eq:c-rich-u-small}, we have
\begin{equation}\label{eq:c-rich-02}
\re \overline{b}\,y(t_1) \geq \re \overline{b}\,b -
|y(t_1)-b|>1-\eps/4.
\end{equation}
We claim that for every $\gamma \in S_{\K}(-a/b,1)$, we have
\begin{equation}\label{eq:c-rich-01}
|a + \gamma\, b| = \left|\gamma - \left(-\frac{a}{b}
\right)\right|=1
\end{equation}
and
\begin{equation}\label{eq:c-rich-03}
\|x+\gamma \,b\,h\|\leq 1 + \eps/4.
\end{equation}
Indeed, the first condition is clear. Let us prove the second one.
If $t\notin U$, then
$$
|x(t)+\gamma\,b\,h(t)|=|x(t)|\leq 1.
$$
If $t\in U$, then
\begin{align*}
|x(t)+\gamma\, b\,h(t)| &\leq |x(t)- a| + |a + \gamma\,b\, h(t)|
\\ &\leq \eps/4 + |a + \gamma\, b\, h(t)|.
\end{align*}
Since $h(t) \in [0,1]$, the number $a + \gamma\, y(t_0)\, h(t)$ is a
convex combination of $a$ and $a + \gamma\, y(t_0)$, so $|a + \gamma
\,y(t_0)\, h(t)| \leq 1$ by Eq.~\eqref{eq:c-rich-01}.

Now, since $0\in \text{co} \left(S_{\K}(-a/b,1)\right)$, we may find
$\gamma_1,\gamma_2\in S_{\K}(-a/b,1)$ and $\lambda\in[0,1]$ such
that $0=\lambda\, \gamma_1+(1-\lambda)\,\gamma_2$. We consider
$$
y^*=\frac{\overline{b}\,\delta_{t_1}|_X}{\left\|\delta_{t_1}|_X\right\|}\in
S_{X^*} \qquad \text{and} \qquad x_i=\frac{x+\gamma_i\,b\,{\tilde
h}}{1+\eps}\in X \quad (i=1,2),
$$
and we observe that Equations \eqref{eq:tildeh} and
\eqref{eq:c-rich-03} give that $x_1,x_2\in B_X$.

Finally, $y\in S(B_X,y^*,2\eps)$ by Eq.~\eqref{eq:c-rich-02},
$$
\left\|x -\left(\lambda\,x_1 + (1-\lambda)\,x_2\right)\right\|
=\frac{\eps}{1+\eps}<\eps,
$$
and, for $i=1,2$,
\begin{align*}
(1+\eps)|y^*(x_i)|&\geq \left|x(t_1) + \gamma_i\, b\,{\tilde
h}(t_1)\right|  \geq |x(t_1) + \gamma_i\, b| - 2\|h-\tilde h\|
\\ & \geq |a + \gamma_i\, b| - |x(t_1)-a| - 2\|h-\tilde
h\|  > 1 - \eps/4 - \eps/2>1-\eps,
\end{align*}
where we have used Equations \eqref{eq:c-rich-u-small},
\eqref{eq:tildeh}, and \eqref{eq:c-rich-01}. Therefore, $x_1,x_2 \in
S(B_X,y^*,2\eps)$.
\end{proof}

When $K$ is perfect, our definition of C-richness coincides with the
definition of richness given in \cite{Ka-Po} and thus every
finite-codimensional subspace of $C(K)$ is C-rich (see
\cite[Proposition~1.2]{Ka-Po}). This is not always the case when $K$
has isolated points. Actually, the following result characterizes
C-rich finite-codimensional subspaces of $C(K)$. We recall that the
\emph{support} of an element $f\in C(K)^*$ (represented by the
regular measure $\mu_f$) is
$$
\text{supp}(f)=\bigcap\left\{C\subset K\ : \ C\ \text{closed},\
|\mu_f|(K\setminus C)=0\right\}.
$$

\begin{prop}\label{prop:description-C-rich}
Let $K$ be a compact Hausdorff   space and let
$f_1,\ldots,f_n\in C(K)^*$. The subspace
$$
Y=\bigcap_{i=1}^n  \ker f_i
$$
is C-rich if and only if $\ \bigcup_{i=1}^n \text{\rm supp}(f_i)$ does
not intersect the set of isolated points of $K$.
\end{prop}

\begin{proof}
We fix a nonempty open subset $U$ of $K$ and $\eps>0$. If
$\bigcup_{i=1}^n \text{supp}(f_i)$ does not contain any isolated
point of $K$, we may consider two cases. Case 1: $U$ contains an
isolated point of $K$ (say, $\tau$). Then $h=\chi_{\{\tau\}} \in
S_{C(K)}$ is a positive $U$-supported function which lies in $Y$, so
$\dist(h,Y)=0<\eps$. Case 2: $U$ does not contain isolated points of
$K.$ In this case one can find a sequence of disjoint open subsets
$U_n \subset U$ and a sequence of positive $h_n \in S_{C(K)}$ with
$\text{supp}(h_n) \subset U_n$. Denote by $q: X \to X/Y$ the natural
quotient map. Since $(h_n)$ tends weakly to $0$ as $n \to \infty$,
$(q(h_n))$ tends weakly to 0. But $X/Y$ is finite-dimensional, so
$\|q(h_n)\|=\dist(h_n,Y) \to 0$ as well, and we can select
$n\in\N$ with $\dist(h_n,Y) < \eps$ and $\text{supp}(h_n)\subseteq
U_n\subset U$. So in both cases $Y$ is C-rich.

Conversely, suppose, for the sake of simplicity, that
$\text{supp}(f_1)$ contains an isolated point $t_0\in K$. Then,
$\mu_{f_1}(\{t_0\})\neq 0$ (if not, $t_0\notin \text{supp}(f_1)$),
and for every positive norm-one function $h$ with support inside
$U$, one has $\dist(h,\ker f_1)\geq |\mu_{f_1}(\{t_0\})|>0$. Hence
$\ker f_1$ is not a C-rich subspace, and neither is~$Y$.
\end{proof}

\section{The counterexamples}\label{sec:counterexample}

Let us recall that $c$ denotes the Banach space of all convergent
scalar sequences $x=(x(1),x(2), \ldots)$ equipped with the
sup-norm. Evidently, $c$ is isometric to $C(K)$ where $K = \mathbb
N \cup \{\infty\}$ is the one-point compactification of $\mathbb
N$. We are now ready for the main result of the paper.

\begin{example}\label{example:main}
{\slshape There exists a Banach space $X$ such that $n(X)=1$ and
$n(X^*)<1$.\ }\newline Indeed, we consider
$$
X=\left\{(x,y,z)\in c\oplus_\infty c \oplus_\infty c\ : \ \lim x +
\lim y + \lim z =0\right\},
$$
which is a C-rich subspace of $c\oplus_\infty c \oplus_\infty c$ by
Proposition~\ref{prop:description-C-rich} and, therefore,
Theorem~\ref{th:C-rich=>fat} gives us that $n(X)=1$. Let us prove
that $n(X^*)<1$. We consider the closed subspace of $X$ given by
$$
Y=\left\{(x,y,z)\in c\oplus_\infty c \oplus_\infty c\ : \ \lim x =
\lim y = \lim z =0\right\}.
$$
Since $Y$ is an $M$-ideal in $c\oplus_\infty c \oplus_\infty c$ (see
\cite[Example~I.1.4(a)]{HWW}), it is a fortiori an $M$-ideal  in $X$ by
\cite[Proposition~I.1.17]{HWW}, meaning that $Y^\perp\equiv
\left(X/Y\right)^*$ is an $L$-summand of $X^*$. Therefore,
$n(X^*)\leq n(Y^\perp)$ by \cite[Proposition~1]{M-P}. But $X/Y$
identifies with the two-dimensional space
$$
\left\{(a,b,c)\in \ell_\infty^{(3)}\ : \ a+b+c=0\right\}
$$
which does not have numerical index~$1$ (in the real case,
Remark~3.6 of \cite{Mc} gives directly the result, since the unit
ball of this space is a hexagon; the complex case follows routinely
from Theorem~3.1 of the same paper).
\end{example}

\begin{remark}
In \cite[Lemma~4.8]{MM-MJ-PA-RP} the reader may find a result which
could be considered as contradictory with the above example. Let us
recall that there is a concept of numerical range for elements of
unital Banach algebras (see \cite[Chapter~1]{B-D1}, for instance).
Given a Banach algebra $A$ with unit $u$, we define the
\emph{algebra numerical range} of an element $a\in A$ by
$$
V(A,a)=\{\varphi(a)\ : \ \varphi\in A^*,\
\|\varphi\|=\varphi(u)=1\}.
$$
We have then a corresponding \emph{algebra numerical radius}
$v(A,a)$ and the corresponding \emph{algebra numerical index}
$n_\text{a}(A)$ of $A$. Given a Banach space $X$, if we consider the
unital Banach algebra $A=L(X)$, it is well-known that
$$
V(L(X),T)=\ecc V(T)
$$
for every $T\in L(X)$ \cite[Theorem~9.4]{B-D1} and thus,
$n(X)=n_\text{a}(L(X))$. It follows from
\cite[Lemma~4.8]{MM-MJ-PA-RP} that
$$
n_\text{a}(L(X))=n_\text{a}(L(X)^{**})
$$
but, in general, $L(X^{**})$ does not coincide with $L(X)^{**}$.
\end{remark}

With just a little bit of work, Example~\ref{example:main} can be
pushed to produce even better counterexamples.

\begin{examples}\label{examples:better} $ $
\begin{enumerate}
\item[(a)] {\slshape There exists a real Banach space $X$
such that $n(X)=1$ and $n(X^*)=0$.\ } Indeed, for every integer
$n\geq 2$, we denote by $Z_n$ the $2$-dimensional real normed space
whose unit ball is the convex hull of the $(2n)^{rm th}$ roots  of unity
(i.e., its unit ball is a regular $2n$-polygon such that one of its
vertices is $(1,0)$). We observe that $Z_n$ is (isometric to) a
subspace of $\ell_\infty^{(n)}$, and it is straightforward,
following the lines of Example~\ref{example:main}, to construct a
C-rich subspace $X_n$ of $c\oplus_\infty c \oplus_\infty \cdots
\oplus_\infty c$, and an $M$-ideal $Y_n$ of $X_n$ such that $X_n/Y_n$
is isometric to $Z_n$. It follows that $n(X_n)=1$ and $n(X_n^*)\leq
n(Z_n)$. Finally, we consider
$$
X:=\biggl[\bigoplus_{n\geq 2}\ X_n\biggr]_{c_0},
$$
and we observe that $n(X)=1$ and $n(X^*)\leq n(X_n)$ for every
$n\geq 2$. But this implies $n(X^*)=0$ since $n(X_n)\longrightarrow
0$ by \cite[Theorem~5]{MarMer}.
\item[(b)] {\slshape There exists a complex Banach space $X$
such that $n(X)=1$ and $n(X^*)=1/\e$.\ } Let $Z$ be a
two-dimensional complex normed space with numerical index $1/\e$
(see \cite[Lemma~32.2]{B-D2}). Then, we may find a family $\{Z_n\}$
of two-dimensional subspaces of $\ell_\infty^{(n)}$ such that the
distance form $Z_n$ to $Z$ goes to $0$. Now, we follow the lines of
the above example to get a Banach space $X$ such that $n(X)=1$ and
$n(X^*)\leq n(Z_n)$ for every $n\in \N$. But the numerical index is
continuous with respect to the distance between Banach spaces
\cite[Proposition~2]{F-M-P}, and so $n(X^*)\leq 1/\e$.
\end{enumerate}
\end{examples}

As we have already mentioned at the beginning of
Section~\ref{sec:C-rich}, the main examples of Banach spaces with
numerical index~$1$ are known to be almost-CL-spaces. Actually, it
is proved in \cite{MarRNP} that every Banach space with numerical
index~$1$ and the RNP is an almost-CL-space and it satisfies that
\begin{equation*}
|x^{**}(x^*)|=1 \qquad \bigl(x^{**}\in \ext[X^{**}],\ x^*\in
\ext[X^*]\bigr).
\end{equation*}
The Example~\ref{example:main} shows that these implications are not
true in general, even for Asplund spaces, as the following result
details. Recall that a \emph{boundary} of $B_{X^*}$ is a subset $C$
of $B_{X^*}$ such that
$$
\|x\|=\max\{\re f(x) \ : \ f\in C\}
$$
for every $x\in X$. The classical boundary of $B_{X^*}$ is the set
$\ext[X^*]$ (consequence of the Hahn-Banach and Krein-Milman
Theorems).

\begin{example}\label{example:Asplund}
{\slshape Let us consider the Banach space given in
Example~\ref{example:main}, i.e.,
$$
X=\left\{(x,y,z)\in c\oplus_\infty c \oplus_\infty c\ : \ \lim x +
\lim y + \lim z =0\right\}.
$$
Then $X$ is an Asplund space, it is lush (and so $n(X)=1$), but the
following properties hold.\ }
\begin{enumerate}
\item[(a)] {\slshape For every boundary $C\subset S_{X^*}$ of $B_{X^*}$,
there exists $x^*\in C$ and $x^{**}\in \ext[X^{**}]$ such that
$|x^{**}(x^*)|<1$.\ } Suppose  that, on the contrary, $C$ is a boundary of
$B_X$ such that $|x^{**}(x^*)|=1$ for every $x^*\in C$ and every
$x^{**}\in \ext[X^{**}]$. Now that  $X$ is a space which does not contain
$\ell_1$, $B_{X^*}$ is the norm-closed convex hull of $C$
\cite[Theorem~III.1]{God}. Therefore, given $T\in L(X^*)$ and
$\eps>0$, we may find $x^*\in C$ and $x^{**}\in \ext[X^{**}]$ such
that
$$
|x^{**}(Tx^*)|=\|Tx^*\|>\|T\|-\eps.
$$
This result together with the fact that $|x^{**}(x^*)|=1$ gives
that $v(T)>\|T\|-\eps$; thus  $n(X^*)=1$, a contradiction.
\item[(b)] In particular, {\slshape there are $x^*\in \ext[X^*]$ and $x^{**}\in
\ext[X^{**}]$ such that $|x^{**}(x^*)|<1$.}
\item[(c)] {\slshape $X$ is not an almost-CL-space.\
} This follows from (a) and \cite[Lemma~3]{MartPaya-CL}.
\end{enumerate}
\end{example}

\begin{remark}
Let us observe that every Asplund space with numerical index~$1$ (in
particular the above example) satisfies the following property:
there is a subset $A$ of $S_{X^*}$ such that $B_{X^*}=\ecc^{w^*}(A)$
and
\begin{equation*}
|x^{**}(x^*)|=1 \qquad \bigl(x^{**}\in \ext[X^{**}],\ x^*\in
A\bigr).
\end{equation*}
This is a consequence of \cite[Lemma~1]{L-M-P}, where $A$ is the set
of all $w^*$-denting points of $B_{X^*}$. The above property is
clearly sufficient for an arbitrary Banach space to have numerical
index~$1$ (see \cite[\S~1]{Mar-ADP}, for instance), but we do not
know if it is also necessary without the Asplundness assumption.
\end{remark}

Once we know that the numerical index of a Banach space and the one
of its dual do not coincide, another natural question could be if
two isometric preduals of a given Banach space should have the same
numerical index. The answer is again negative as the following
result shows.

\begin{example}
{\slshape There is a Banach space $Z$ with two isometric preduals
$X_1$ and $X_2$ such that $n(X_1)$ and $n(X_2)$ are not equal.\ }
Indeed, let
$$
X_1=\left\{(x,y,z)\in c\oplus_\infty c \oplus_\infty c\ : \ \lim x +
\lim y + \lim z =0\right\}
$$
and
$$
X_2=\left\{(x,y,z)\in c\oplus_\infty c \oplus_\infty c\ : \ x(1) +
y(1) + z(1) =0\right\}.
$$
By Example~\ref{example:main}, $n(X_1)=1$. Since the two dimensional
space
$$
\left\{(a,b,c)\in \ell_\infty^{(3)}\ : \ a+b+c=0\right\}
$$
is isometric to an $M$-summand of $X_2$, it follows that
$n(X_2)<1$ (see \cite[Proposition~1]{M-P} and
\cite[Theorem~3.1]{Mc}). Finally, the fact that $X_1^*$ and $X_2^*$
are isometric is straightforward.
\end{example}

\section{A positive result}\label{sec:positive}
As a straightforward application of the inequality
\eqref{eq:inequality}, i.e., $n(X^*)\leq n(X)$, it is clear that
$n(X)=n(X^*)$ for every reflexive space $X$. This equality also
holds when $X$ is a Banach space such that $n(X^*)=1$, in particular
when $X$ is an $L$- or an $M$-space. Besides  these elementary
results, it is also true that $n(X)=n(X^*)$ when $X$ is a
$C^*$-algebra or a von Neumann predual (see \cite{Hur} and
\cite[pp.~202]{K-M-RP}).

We finish the paper by showing another particular case where
inequality~\eqref{eq:inequality} becomes an equality.

\begin{prop}
Let $X$ be a Banach space with the RNP. If $n(X)=1$, then
$n(X^*)=1$.
\end{prop}

\begin{proof}
By \cite[Lemma~1]{L-M-P}, we have that $|x^*(x)|=1$ for every
extreme point $x^*$ of $B_{X^*}$ and every denting point $x\in B_X$.
Therefore, \cite[Proposition~2.1]{Dutta-Rao} (or
\cite[Proposition~3.5]{Sharir}) gives us that
\begin{equation}\label{eq:RNP-1}
|x^{***}(x)|=1
\end{equation}
for every $x^{***}\in \ext[X^{***}]$ and every denting point $x\in
B_X$. Now, we fix $T\in L(X^*)$ and $\eps>0$. Since $X$ has the RNP,
$B_{X^{**}}$ is the weak$^*$-closed convex hull of the set of
denting points of $B_X$, and we may find  a denting point $x$ such
that
$$
\|T^* x\|>\|T\| - \eps.
$$
Then, we may find $x^{***}\in \ext[X^{***}]$ such that
$$
|x^{***}(T^* x)|=\|T^* x \|> \|T\| - \eps.
$$
This fact, together with Eq.~\eqref{eq:RNP-1}, implies that
$\|T^*\|-\eps\leq v(T^*)$. By  letting $\eps\downarrow 0$, we
have
\begin{equation*}
\|T\|=\|T^*\|=v(T^*)=v(T).\qedhere
\end{equation*}
\end{proof}

We do not know if $n(X)=n(X^*)$ for every Banach space with the
Radon-Nikod\'{y}m property.


\begin{thebibliography}{99}

\bibitem{Bauer} \textsc{F.~L.~Bauer},
On the field of values subordinate to a norm, \emph{Numer. Math.}
\textbf{4} (1962), 103--111.

\bibitem{B-K} \textsc{H.~F.~Bohnenblust and S.~Karlin},
Geometrical properties of the unit sphere in Banach algebras,
\emph{Ann. of Math.}, \textbf{62} (1955), 217--229.

\bibitem{B-D1} \textsc{F.~F.~Bonsall and J.~Duncan},
\emph{Numerical Ranges of Operators on Normed Spaces and of Elements
of Normed Algebras}, London Math. Soc. Lecture Note Series
\textbf{2}, Cambridge, 1971.

\bibitem{B-D2} \textsc{F.~F.~Bonsall and J.~Duncan},
\emph{Numerical Ranges II}, London Math. Soc. Lecture Note Series
\textbf{10}, Cambridge, 1973.


\bibitem{C-D-Mc} \textsc{M.~J.~Crabb, J.~Duncan, and C.~M.~McGregor},
Mapping theorems and the numerical radius, \emph{Proc. London Math.
Soc.} \textbf{25} (1972), 486--502.


\bibitem{D-Mc-P-W} \textsc{J.~Duncan, C.~McGregor, J.~Pryce, and
A.~White}, The numerical index of a normed space, \emph{J. London
Math. Soc.} \textbf{2} (1970), 481--488.

\bibitem{Dutta-Rao} \textsc{S.~Dutta and T.~S.~S.~R.~K.~Rao}, On
weak$^*$-extreme points in Banach spaces, \emph{J. Convex Anal.}
\textbf{10} (2003), 531--539.

\bibitem{Ed-Dari} \textsc{E.~Ed-Dari},
On the numerical index of Banach spaces, \emph{Linear Algebra Appl.}
403 (2005), 86--96.

\bibitem{Ed-Dari-Khamsi} \textsc{E.~Ed-Dari and M.~A.~Khamsi},
The numerical index of the $L_p$ space, \emph{Proc. Amer. Math.
Soc.}, to appear.

\bibitem{F-M-P} \textsc{C.~Finet, M.~Mart\'{\i}n, and R.~Pay\'{a}},
Numerical index and renorming, \emph{Proc. Amer. Math. Soc.}
\textbf{131} (2003), 871--877.

\bibitem{Full} \textsc{R.~E.~Fullerton},
\emph{Geometrical characterization of certain function spaces}. In:
Proc. Inter. Sympos. Linear spaces (Jerusalem 1960), pp. 227--236.
Pergamon, Oxford 1961.

\bibitem{Gli} \textsc{B.~W.~Glickfeld},
On an inequality of Banach algebra geometry and semi-inner-product
space theory, \emph{Illinois J. Math.}, \textbf{14} (1970), 76--81.

\bibitem{God} \textsc{G.~Godefroy},
Boundaries of a convex set an interpolation sets, \emph{Math. Ann.}
\textbf{277} (1987), 173--184.


\bibitem{Hal} \textsc{P.~Halmos},
\emph{A Hilbert space problem book}, Van Nostrand, New York, 1967.

\bibitem{HWW} \textsc{P.~Harmand, D.~Werner, and D.~Werner},
\emph{$M$-ideals in Banach spaces and Banach algebras}, Lecture
Notes in Math. \textbf{1547}, Springer-Verlag, Berlin 1993.

\bibitem{Hur} \textsc{T.~Huruya}, The normed space numerical index
of $C^*$-algebras, \emph{Proc. Amer. Math. Soc.} \textbf{63} (1977),
289--290.

\bibitem{Ka-Po} \textsc{V.~M.~Kadets and M.~M.~Popov}, The Daugavet property
for narrow operators in rich subspaces of $C[0,1]$ and $L_1[0,1]$,
\emph{St. Petersburg Math. J.} \textbf{8} (1997), 571--584.

\bibitem{K-M-RP} \textsc{A.~Kaidi, A.~Morales, and
A.~Rodr\'{\i}guez-Palacios}, Geometrical properties of the product of a
$C^*$-algebra, \emph{Rocky Mountain J. Math.} \textbf{31} (2001),
197--213.

\bibitem{Lima2} \textsc{\AA.~Lima},
\emph{Intersection properties of balls in spaces of compact
operators}, Ann. Inst. Fourier, Grenoble \textbf{28} (1978), 35--65.



\bibitem{L-M-P} \textsc{G.~L\'{o}pez, M.~Mart\'{\i}n, and R.~Pay\'{a}},
Real Banach spaces with numerical index 1, \emph{Bull. London Math.
Soc.} \textbf{31} (1999), 207--212.

\bibitem{Lumer} \textsc{G.~Lumer},
Semi-inner-product spaces, \emph{Trans. Amer. Math. Soc.}
\textbf{100} (1961), 29--43.

\bibitem{Mar} \textsc{M.~Mart\'{\i}n},
A survey on the numerical index of a Banach space, \emph{Extracta
Math.} \textbf{15} (2000), 265--276.


\bibitem{MarRNP} \textsc{M.~Mart\'{\i}n},
Banach spaces having the Radon-Nikod\'{y}m property and numerical index
$1$, \emph{Proc. Amer. Math. Soc.} \textbf{131} (2003), 3407--3410.

\bibitem{Mar-ADP} \textsc{M.~Mart\'{\i}n},
The alternative Daugavet property for $C^*$-algebras and
$JB^*$-triples, preprint. Available at
\texttt{http://arXiv.org/abs/math.FA/0411555}

\bibitem{MarMer} \textsc{M.~Mart\'{\i}n and J.~Mer\'{\i}}, Numerical index
of some polyhedral norms on the plane, \emph{preprint}.

\bibitem{MaMeRo} \textsc{M.~Mart\'{\i}n, J.~Mer\'{\i}, and A.~Rodr\'{\i}guez-Palacios},
Finite-dimensional Banach spaces with numerical index zero,
\emph{Indiana U. Math. J.} \textbf{53} (2004), 1279--1289.

\bibitem{MaOi} \textsc{M.~Mart\'{\i}n and T.~Oikhberg},
An alternative Daugavet property, \emph{J. Math. Anal. Appl.}
\textbf{294} (2004), 158--180.

\bibitem{M-P} \textsc{M.~Mart\'{\i}n and R.~Pay\'{a}},
Numerical index of vector-valued function spaces, \emph{Studia
Math.} \textbf{142} (2000), 269--280.

\bibitem{MartPaya-CL} \textsc{M.~Mart\'{\i}n and R.~Pay\'{a}},
On CL-spaces and almost-CL-spaces, \emph{Ark. Mat.} \textbf{42}
(2004), 107--118.

\bibitem{M-V} \textsc{M.~Mart\'{\i}n and A.~R.~Villena},
Numerical index and Daugavet property for $L_\infty(\mu,X)$,
\emph{Proc. Edinburgh Math. Soc.} \textbf{46} (2003), 415--420.

\bibitem{MM-MJ-PA-RP} \textsc{J.~Mart\'{\i}nez-Moreno, J.~F.~Mena-Jurado,
R.~Pay\'{a}-Albert, and A.~Rodr\'{\i}guez-Palacios}, An approach to numerical
ranges without Banach algebra theory, \emph{Illinois J. Math.}
\textbf{29} (1985), 609--626.

\bibitem{Mc} \textsc{C.~M.~McGregor},
Finite dimensional normed linear spaces with numerical index $1$,
\emph{J. London Math. Soc.} \textbf{3} (1971), 717--721.

\bibitem{Oik04} \textsc{T.~Oikhberg},
Spaces of operators, the $\psi$-Daugavet property, and numerical
indices, \emph{Positivity} (to appear).


\bibitem{Reis} \textsc{S.~Reisner},
Certain Banach spaces associated with graphs and CL-spaces with
$1$-unconditonal bases, \emph{J. London Math. Soc.} \textbf{43}
(1991), 137--148.

\bibitem{Sharir} \textsc{M.~Sharir}, Extremal structure in operator
spaces, \emph{Trans. Amer. Math. Soc.} \textbf{186} (1973), 91--111.

\bibitem{Toe} \textsc{O.~Toeplitz},
Das algebraische Analogon zu einem Satze von Fejer, \emph{Math. Z.}
\textbf{2} (1918), 187--197.

\bibitem{WerJFA} \textsc{D.~Werner},
The Daugavet equation for operators on function spaces, \emph{J.
Funct. Anal.} \textbf{143} (1997), 117--128.

\end{thebibliography}
\end{document}